\documentstyle[12pt]{article}
\begin{titlepage}
\title{
A theta function for hyperbolic surfaces with cusps}

\author{Ulrich Bunke\thanks{Humboldt-Universit\"at zu Berlin, Institut f\"ur
Reine Mathematik (SFB288), Ziegelstr. 13a, Berlin 10099.
E-mail:ubunke@mathematik.hu-berlin.de
       } and
Martin
Olbrich\thanks{Humboldt-Universit\"at zu Berlin, Institut f\"ur Reine
Mathematik (SFB288), Ziegelstr. 13a, Berlin 10099.
E-mail:olbrich@mathematik.hu-berlin.de  }
}
\end{titlepage}

\begin{document}

\newcommand{\R}{{\bf R}}
\newcommand{\Z}{{\bf Z}}
\newcommand{\C}{{\bf C}}
\newcommand{\N}{{\bf N}}
\maketitle
\newtheorem{prop}{Proposition}[section]
\newtheorem{lem}[prop]{Lemma}
\newtheorem{ddd}[prop]{Definition}
\newtheorem{theorem}[prop]{Theorem}
\newtheorem{kor}[prop]{Corollary}
\newtheorem{ass}[prop]{Assumption}
\newtheorem{con}[prop]{Conjecture}

\begin{abstract}
For a Riemann surface with cusps we define a theta function using the
eigenvalues of the
Laplacian and the singularities of the scattering determinant. We provide its
meromorphic
continuation and discuss its singularities.
\end{abstract}

\tableofcontents
\section{Theta functions}

Let $M$ be a complete Riemann surface of constant negative curvature $-1$
of finite volume. Consider the unique selfadjoint extension
$\Delta_M$ of the Laplacian and form $A:=\sqrt{\Delta_M-\frac{1}{4}}$,
where we take the square root with positive imaginary part.
Let $P_d$ be the projection onto the subspace of $L^2(M)$ spanned
by the $L^2$-eigenfunctions of $\Delta_M$, i.e. onto the discrete
subspace. The (eigenvalue) theta function of $M$ is defined by
\begin{ddd}
$$\theta_d(t):=Tr P_d e^{\imath t A}P_d ,\quad Im(t)>0.$$
\end{ddd}
It is known, that the number of eigenvalues of $\Delta_M$
(counted with multiplicity)
being smaller than $\lambda$ is bounded by $C\lambda^2$ for some
$C<\infty$. Hence, $\theta_d(t)$ is a well defined holomorphic
function on the upper half plane.

If $M$ is compact, this $\theta$-function was introduced by Cartier/Voros
\cite{cartiervoros90}. They provide a meromorphic extension
of $\theta_d(t)$ to the
complex plane and give a complete description of its singularities.
Their method starts with the Selberg Zeta-function of $M$.
The theta-function is derived by a certain contour integral
of the logarithmic derivative of the Selberg zeta-function.

In Bunke/Olbrich/Juhl \cite{bunkeolbrich93} we developed
another method to obtain the meromorphic extension starting
with a distributional trace formula. In the present paper
we use the trace formula in a similar way.
If $M$ is non-compact, we are not able to provide an
extension of $\theta_d(t)$ alone. The point is, that the eigenvalue
theta-function enters into the trace formula in combination
with another (scattering) theta-function obtained from
the singularities of the scattering determinant.

Let $S(z)$ be the scattering determinant of $M$. $S(z)$
is a meromorphic function of finite order, satisfying
$$ S(z)=S(-z)^{-1},\quad \bar{S(z)}=S(\bar{z})^{-1},\quad z\in\C.$$
Let $\Sigma$ be the set of singularities $\sigma$ of $S(z)$ with
$Re(\sigma)\ge 0$, $Im(\sigma)>0$. Let $m_\sigma$ ($-m_\sigma,
2m_\sigma,-2m_\sigma$) be the order of the pole with $Re(\sigma)>0$
(zero with $Re(\sigma)>0$, pole with $Re(\sigma)=0$, zero with $Re(\sigma)=0$).
It is known, that for some $C<\infty$ and all $\sigma\in\Sigma$ $Im(\sigma)<C$.
Moreover, there are at most finitely many zeros on $\Sigma$ and
the total multiplicity of all singularities with $Re(\sigma)<\lambda$
is bounded by $C\lambda^2$ for some $C<\infty$. Thus,
\begin{ddd}
$$\theta_s(t):=\sum_{\sigma\in\Sigma} m_\sigma e^{\imath t\sigma},\quad
Im(t)>0$$
\end{ddd}
is a well defined holomorphic function on the upper half plane.

The main object of the present paper is the theta-function associated
to $M$:
\begin{ddd}
$$\theta(t):=\theta_d(t)+\theta_s(M),$$
\end{ddd}
which is initially defined on the upper half-plane.

In the  second section we state the version of the Selberg trace
formula suitable for our purpose. In the
third section we give a meromorphic extension of $\theta(t)$ to the
lower half-plane across $(0,\infty)$. Here, the trace formula
can be applied immediately. In the forth section  we
discuss the necessary modifications needed to provide a
meromorphic extension of the theta-function across the negative
real axis. Finally, in the fifth section we collect all results
together and state our main theorem about the meromorphic extension
of the (modified) theta-function to the complex plane, and we
give a complete description of its singularities.

\section{The Selberg trace formula}

Let $\hat{\phi}$ be the Fourier transform
$$\hat{\phi}(\lambda)=\int_{-\infty}^\infty e^{\imath\lambda t} \phi(t) dt$$
of an even function $\phi$ satisfying for some $\delta>0$, $C<\infty$
\begin{eqnarray*}
|\phi(t)|&<&Ce^{(-1/2+\delta)|t|}\\
|\frac{d\hat{\phi}}{d\lambda}|&<&\frac{C}{(1+|\lambda|)^{3+\delta}}.
\end{eqnarray*}
We introduce the hyperbolic contribution.
\begin{ddd}
$$c_h(\phi)=\sum_{c \mbox{\scriptsize -closed
geodesic}}\frac{l_c}{2n_csinh(l_c/2)}\phi(l_c),$$
where $l_c$ is the length of the geodesic $c$ and $n_c$ is its multiplicity.
\end{ddd}
The contribution of the identity is given by
\begin{ddd}
$$c_e(\phi)=-\frac{vol(M)}{4\pi}\int_{-\infty}^\infty \frac{\frac{d}{dt}
\phi(t)}{sinh(t/2)}dt.$$
\end{ddd}
Note that the integral is well defined at zero, since $\phi$ is even.
The Selberg trace formula, as it can be found e.g in Lax/Phillips
\cite{laxphillips76},
reads
\begin{eqnarray*}
\sum_{\lambda \mbox{\scriptsize\it eigenvalue of $A$}}
\hat{\phi}(\lambda)&=&c_e(\phi)+c_h(\phi)\\
&-&r\:ln(2) \phi(0) + r/2 \hat{\phi}(0)\\
&-&\frac{r}{2\pi} \int_{-\infty}^\infty
\hat{\phi}(\lambda)\frac{\dot{\Gamma}
(1+\imath\lambda)}{\Gamma(1+\imath\lambda)}d\lambda\\
&-&\frac{\imath}{4\pi}\int_{-\infty}^\infty
\hat{\phi}(\lambda)
\frac{\dot{S}(\lambda)}{S(\lambda)} d\lambda.
\end{eqnarray*}
Set $\phi(t)=\psi(t)+\psi(-t)$ for some $\psi\in C_c^\infty(0,\infty)$.
We can write all terms of the trace formula as applications of distributions
to $\psi$.

Let
\begin{ddd}
$$\tilde{\theta}_d(t):=Tr P_d e^{-\imath t A} P_d,\quad Im(t)<0.$$
\end{ddd}
$\tilde{\theta}_d(t)$ is a holomorphic function on the lower half-plane.
The left-hand side of the trace formula can be written as
$$<\theta_d(t+\imath 0)+\tilde{\theta}_d(t-\imath 0),\psi>.$$
The hyperbolic contribution is
$$<\sum_{c\mbox{\scriptsize\it -closed
geodesic}}\frac{l_c}{2n_csinh(l_c/2)}\delta(t-l_c),\psi>$$
and the contribution of the identity is
$$<- \frac{vol(M)}{4\pi}\frac{cosh(t/2)}{sinh^2(t/2)},\psi>.$$
The last two terms of the trace formula are less trivial.
We consider first
$$-\frac{r}{2\pi} \int_{-\infty}^\infty
\hat{\psi}(\lambda)\frac{\dot{\Gamma}
(1+\imath\lambda)}{\Gamma(1+
\imath\lambda)}d\lambda.$$
Since $\Gamma(1+\imath\lambda)$
is of finite order and
$\hat{\psi}$ vanishes exponentially in the upper half-plane, we can apply
the Cauchy integral formula.
$\Gamma(1+\imath\lambda)$ has poles of first order in the points
$\lambda=\imath k$, $k=1,2,\dots$.
Thus
\begin{eqnarray*}
-\frac{r}{2\pi} \int_{-\infty}^\infty
\hat{\psi}(\lambda)\frac{\dot{\Gamma}(1+\imath
\lambda)}{\Gamma(1+\imath\lambda)}d\lambda
&=&r \sum_{k=1}^\infty
\hat{\psi}(\imath k)\\
&=&<\frac{re^{-t}}{1-e^{-t}},\psi>.
\end{eqnarray*}
In an analogous manner
$$-\frac{r}{2\pi} \int_{-\infty}^\infty
\hat{\psi_-}(\lambda)\frac{\dot{\Gamma}
(1+\imath\lambda)}{\Gamma(1+\imath\lambda)}d\lambda=0,$$
where $\psi_-(t):=\psi(-t)$.
The term involving the scattering determinant is even more complicate.
Fortunately, the following equation was shown by M\"uller \cite{mueller912}:
$$
\frac{\imath}{4\pi}\int_{-\infty}^\infty \hat{\phi}(\lambda)
\frac{\dot{S}(\lambda)}{S(\lambda)} d\lambda\\
=\sum_{\sigma\in\Sigma} m_\sigma (\hat{\psi}(\sigma)+\hat{\psi}(-\bar{\sigma}))
$$
Let
\begin{ddd}
$$\tilde{\theta}_s(t):=\bar{\theta}_s(\bar{t}),\quad Im(t)<0.$$
\end{ddd}
$\tilde{\theta}_s(t)$ is a holomorphic function on the lower half-plane.
We can write
$$\frac{\imath}{4\pi}\int_{-\infty}^\infty \hat{\phi}(\lambda)
\frac{\dot{S}(\lambda)}{S(\lambda)} d\lambda=<\theta_s(t+\imath
0)+\tilde{\theta}_s(t-\imath 0),\psi>.$$
Let $\tilde{\theta}:=\tilde{\theta}_d+\tilde{\theta}_s$.
The following Lemma is now obvious.
\begin{lem}\label{l1}
The following equation of distributions holds on $(0,\infty)$
\begin{eqnarray*}\theta(t+\imath 0)+\tilde{\theta}(t-\imath 0)&=&-\frac{
vol(M)}{4\pi}\frac{cosh(t/2)}{sinh^2(t/2)}\\
 &+& \sum_{c\mbox{\scriptsize -closed
geodesic}}\frac{l_c}{2n_csinh(l_c/2)}\delta(t-l_c)\\
 &+&\frac{r}{1-e^{-t}}.\end{eqnarray*}
\end{lem}

\section{Analytic continuation across $(0,\infty)$}

\begin{theorem}\label{th3}
$\theta(t)$ has a meromorphic continuation across $(0,\infty)$ to
the lower half-plane
\begin{equation}\label{ee1}\theta(t)=-\tilde{\theta}(t)-
\frac{vol(M)}{4\pi}\frac{cosh(t/2)}{sinh^2(t/2)}+
\frac{r}{1-e^{-t}},\quad Im(t)<0. \end{equation}
The singularities are
\begin{itemize}
\item first order poles at $l_c$, $c$ a closed geodesic, with residue
$\frac{l_c}{4\pi n_csinh(l_c/2)}$ and
\item and second order poles at $-2 \pi\imath k$, $k=1,2,\dots$ with the same
singular part as
      $$-\frac{vol(M)}{4\pi}\frac{cosh(t/2)}{sinh^2(t/2)}+\frac{r}{1-e^{-t}}.$$
\end{itemize}
\end{theorem}
{\it Proof:}
We use (\ref{ee1}) in order to define $\theta(t)$ as a distribution on
$U:=\C\setminus ((-\infty,0]\cup -2 \pi \imath \N)$.
Then we compute $\bar{\partial}\theta$. The interesting region is near
$(0,\infty)$.
For $\phi\in C_c^\infty(U)$ we define
\begin{eqnarray*}<\theta,\phi>&=&
lim_{\epsilon\to+0}\int_\epsilon^\infty
\int_{-\infty}^{\infty}
\theta
(t+\imath u)\phi(t+
\imath u) dt du\\
&-&lim_{\epsilon\to+
0}\int_\epsilon^\infty\int_{
-\infty}^{\infty}
\theta(t-\imath u)\phi(t-
\imath u) dt du.
\end{eqnarray*}
In view of the polynomial bounds on the discrete
spectrum and the singularities
of the scattering determinant this defines a distribution.
We compute, using Lemma \ref{l1},
\begin{eqnarray*}
<\bar{\partial}\theta,\phi>&=&<\theta,\bar{\partial}^\ast\phi>\\
&=&lim_{\epsilon\to+0}\int_\epsilon^\infty\int_{-\infty}^{\infty}
\theta(t+\imath u)\frac{1}{2}\left[-\frac{\partial}{\partial (\imath
u)}-\frac{\partial}{\partial t}\right]\phi(t+\imath u) dt du\\
&+&lim_{\epsilon\to+0}\int_\epsilon^\infty\int_{-\infty}^{\infty}
\theta(t-\imath u)\frac{1}{2}\left[\frac{\partial}{\partial (\imath
u)}-\frac{\partial}{\partial t}\right]\phi(t-\imath u) dt du\\
&=&lim_{\epsilon\to+0}\int_{-\infty}^{\infty} \theta(t+\imath
\epsilon)\frac{1}{2} \phi(t+\imath \epsilon) dt\\
&-&lim_{\epsilon\to+0}\int_{-\infty}^{\infty} \theta(t-\imath
\epsilon)\frac{1}{2} \phi(t-\imath \epsilon) dt\\
&=&\frac{1}{2}<\theta(t+\imath 0)+\tilde{\theta}(t-\imath
0)+\frac{vol(M)}{4\pi}\frac{cosh(t/2)}{sinh^2(t/2)}-\frac{r}{1-e^{-t}},\phi>\\
&=&<\sum_{c\mbox{\scriptsize \it-closed
geodesic}}\frac{l_c}{4n_csinh(l_c/2)}\delta(t-l_c),\phi>.
\end{eqnarray*}
Thus, $\theta(t)$ is holomorphic on $U$ except the points $l_c$, $c$ closed
geodesic.
Since $\bar{\partial}\frac{1}{t}=\pi\delta(t)$, at these points are first order
poles
with the residues claimed in the Theorem. Since $\tilde{\theta}$ is holomorphic
for $Im(t)<0$
the singularities in the lower half plane come from the addittional terms
$$-\frac{vol(M)}{4\pi}\frac{cosh(t/2)}{sinh^2(t/2)}+\frac{r}{1-e^{-t}}.$$
$\Box$\newline
Unfortunately, this method cannot be applied in order to provide an extension
of $\theta$ across the negative real half-axis. We do not know the necessary
trace formula as Lemma \ref{l1}.

\section{$\theta$ near $(-\infty,0)$}

Note that the trace formula \ref{l1} involves the distributions
$$\theta_s(t+\imath 0)+\tilde{\theta}_s(t-\imath 0), \quad\theta_d(t+\imath
0)+\tilde{\theta}_d(t-\imath 0).$$
Obviously
$$\theta_d(t+\imath 0)+\tilde{\theta}_d(t-\imath 0)=\theta_d(-t+\imath
0)+\tilde{\theta}_d(-t-\imath 0).$$
But
$$\theta_s(t+\imath 0)+\tilde{\theta}_s(t-\imath 0)\not=\theta_s(-t+\imath
0)+\tilde{\theta}_s(-t-\imath 0).$$
We will apply the method of Cramer \cite{cramer19}, Jorgenson/Lang
\cite{jorgensonlang93}
in order to provide a meromorphic extension of
$$V(t):=\tilde{\theta}_s(-{t})-\theta_s(t),\quad Im(t)>0.$$
Note that
$$V(t)=\sum_{\sigma\in\Sigma} m_\sigma (e^{\imath t \bar{\sigma}}-e^{\imath t
\sigma}).$$
We write the scattering determinant as $S(\lambda)=G(\lambda)L(\lambda)$, where
$$G(\lambda):=e^{a+b\lambda}\left[\sqrt{\pi}\frac{
\Gamma(\imath
\lambda)}{\Gamma(1/2+\imath\lambda)}
\right]^r.$$
(recall that $r$ is the number of cusps).
By Hejhal, \cite{hejhal83}, ch.8, 3.35, $L(\lambda)$ has a representation
as a certain Dirichlet series of the form
$$L(\lambda):=1+\sum_{q\in Q} \frac{c_q}{q^{\imath \lambda}}.$$
Here, $Q$ is a certain set of positive real numbers $>1$ and $c_q$
are real 'multiplicities'. In order to obtain the $1$ in front of the
sum one has to adjust $a,b\in \R$ appropriately.
The Dirichlet series converges for $Im(\lambda)$ small enough.
Let $R:=\{(n,q_1,\dots,q_n)| n\in\N, q_i\in Q, i=1,\dots,n\}$.
For $p\in R$ define the multiplicity
$$c(p):=\frac{(-1)^{n+1}}{n}\prod_{q\in p} c_q$$
and $$|p|:=\prod_{q\in p} q.$$
Then $ln\:L(\lambda)$ also has a Dirichlet series representation
$$ln\:L(\lambda):=\sum_{p\in R} \frac{c(p)}{|p|^{\imath \lambda}},$$
which converges uniformly and absolutely in any half-plane
$Im(\lambda)\le -A$ for some $\infty>A>0$.
Since $S$ and $G$ are meromorphic functions of finite order, $L$
is also a meromorphic function of finite order.
$L$ satisfies the functional equation
\begin{equation}\label{ee3}L(\lambda)=
\frac{1}{G(-\lambda)G(\lambda)}\frac{1}{L(-\lambda)}.\end{equation}
In Jorgenson/Lang it was shown that $V(t)$ has a representation
as a contour integral. Let $\gamma$ be the piecewise straight path
going from $\imath A+\infty$ to $\imath A$, then down to $-\imath A$ and
then back to $-\imath A+\infty$. The down path avoids the singularities
of $S(z)$ on $Re(\lambda)=0$ on positive real-part side.
By Lang/Jorgenson \cite{jorgensonlang93}, Thm. 2.3, we have
$$2\pi\imath V(t)=\int_\gamma e^{\imath t \lambda}
\frac{\dot{S}(\lambda)}{S(\lambda)}+2\pi\imath \sum_{\sigma\in\Sigma,
Re(\sigma)=0} m_\sigma (e^{-\imath t \sigma}-e^{\imath t \sigma}).$$
We will study the integral further.
Note that
$$\frac{\dot{S}(\lambda)}{S(\lambda)}=
\frac{\dot{G}(\lambda)}{G(\lambda)}+\frac{\dot{L}(\lambda)}{L(\lambda)}$$
$G$ has no singularities in the region encircled by $\gamma$. Thus
$$\int_\gamma e^{\imath t \lambda}
\frac{\dot{S}(\lambda)}{S(\lambda)}=\int_\gamma e^{\imath t \lambda}
\frac{\dot{L}(\lambda)}{L(\lambda)}.$$
Let $h_1$ be the global holomorphic function
$$h_1(t):=-\imath\int_{-A}^Ae^{-t\lambda}\frac{\dot{L}(\imath
\lambda)}{L(\imath\lambda)} d\lambda$$
(avoiding the singularities  as above).
Then
\begin{eqnarray}
\int_\gamma e^{\imath t \lambda}
\frac{\dot{L}(\lambda)}{L(\lambda)}&=&h_1(t)\nonumber\\
&-&e^{-At}\int_0^\infty e^{\imath t\lambda} \frac{\dot{L}(\imath
A+\lambda)}{L(\imath A+\lambda)} d\lambda\label{e2e}\\
&+&e^{At}\int_0^\infty e^{\imath t\lambda} \frac{\dot{L}(-\imath
A+\lambda)}{L(-\imath A+\lambda)} d\lambda\nonumber.
\end{eqnarray}
We employ now the functional equation (\ref{ee3}) in order to rewrite
(\ref{e2e}).
\begin{eqnarray}
-e^{-At}\int_0^\infty e^{\imath t\lambda} \frac{\dot{L}(\imath
A+\lambda)}{L(\imath A+\lambda)} d\lambda
&=&-e^{-At}\int_0^\infty e^{\imath t\lambda} \frac{\dot{L}(-\imath
A-\lambda)}{L(-\imath A-\lambda)} d\lambda\nonumber\\
&-&e^{-At}\int_0^\infty e^{\imath t\lambda} \frac{\dot{G}(-\imath
A-\lambda)}{G(-\imath A-\lambda)} d\lambda\label{ee5}\\
&+&e^{-At}\int_0^\infty e^{\imath t\lambda} \frac{\dot{G}(\imath
A+\lambda)}{G(\imath A+\lambda)} d\lambda. \nonumber
\end{eqnarray}
In order to evaluate the integrals involving $L$ we use the representation of
$ln\:L$ as
a Dirichlet series.
\begin{eqnarray*}
&&-e^{-At}\int_0^\infty e^{\imath t\lambda} \frac{\dot{L}(-\imath
A-\lambda)}{L(-\imath A-\lambda)} d\lambda\\
&=&-e^{-At} ln\: L(-\imath A)-e^{-At}\imath t   \int_0^\infty e^{\imath
t\lambda} ln\:L(-\imath A-\lambda)  d\lambda\\
&=&-e^{-At} ln\: L(-\imath A)-e^{-At}\imath t   \int_0^\infty e^{\imath
t\lambda} \sum_{p\in R} \frac{c(p)}{|p|^{A-\imath \lambda }}  d\lambda\\
&=&-e^{-At} ln\: L(-\imath A)-e^{-At}\imath t \sum_{p\in R} c(p) \int_0^\infty
e^{\imath t\lambda}  e^{(-A+\imath \lambda) ln\:|p|}  d\lambda\\
&=&-e^{-At} ln\: L(-\imath A)+e^{-At} t \sum_{p\in R}
\frac{c(p)}{|p|^{A}(t+ln\:|p|)} .
\end{eqnarray*}
Analogously
\begin{eqnarray*}
&&e^{At}\int_0^\infty e^{\imath t\lambda} \frac{\dot{L}(-\imath
A+\lambda)}{L(-\imath A+\lambda)} d\lambda\\
&=&-e^{At} ln\: L(-\imath A)-e^{At}\imath t   \int_0^\infty e^{\imath t\lambda}
ln\:L(-\imath A+\lambda)  d\lambda\\
&=&-e^{At} ln\: L(-\imath A)-e^{At}\imath t   \int_0^\infty e^{\imath t\lambda}
\sum_{p\in R} \frac{c(p)}{|p|^{A+\imath \lambda }}  d\lambda\\
&=&-e^{At} ln\: L(-\imath A)-e^{At}\imath t \sum_{p\in R} c(p) \int_0^\infty
e^{\imath t\lambda}  e^{(-A-\imath \lambda) ln\:|p|}  d\lambda\\
&=&-e^{At} ln\: L(-\imath A)+e^{At} t \sum_{p\in R}
\frac{c(p)}{|p|^{A}(t-ln\:|p|)} .
\end{eqnarray*}
Define the global holomorphic function
$$h_2(t):=-e^{-At} ln\: L(-\imath A)-e^{At} ln\: L(-\imath A)$$
and the global meromorphic function
\begin{equation}\label{er3}W(t):=e^{-At} t \sum_{p\in R}
\frac{c(p)}{|p|^{A}(t+ln\:|p|)}+e^{At} t \sum_{p\in R}
\frac{c(p)}{|p|^{A}(t-ln\:|p|)}.\end{equation}
Then
$$-e^{-At}\int_0^\infty e^{\imath t\lambda} \frac{\dot{L}(\imath
A+\lambda)}{L(\imath A+\lambda)} d\lambda\\
+e^{At}\int_0^\infty e^{\imath t\lambda} \frac{\dot{L}(-\imath
A+\lambda)}{L(-\imath A+\lambda)} d\lambda\\
=h_2(t)+W(t).
$$
Note that
\begin{equation}\label{ee4}
W(t)=W(-t)\ ,\quad h_2(t)=h_2(-t).
\end{equation}

It remains to study the integrals involving the $G$-factors.
We have
$$\frac{\dot{G}(\lambda)}{G(\lambda)}=b+\imath r\left(
\frac{\dot{\Gamma}(\imath \lambda)}{\Gamma(\imath \lambda)} -
\frac{\dot{\Gamma}(1/2+\imath \lambda)}{\Gamma(1/2+\imath \lambda)}\right).$$
The term involving $b$ will cancel out eventually (see (\ref{ee5})).

We consider
$$
-e^{-At}\int_0^\infty e^{\imath t\lambda} \frac{\dot{G}(-\imath
A-\lambda)}{G(-\imath A-\lambda)} d\lambda.
$$
We evaluate
$$I_1(t):=-\imath e^{-At}\int_0^\infty e^{\imath t\lambda}
\frac{\dot{\Gamma}(A-\imath \lambda)}{\Gamma(A-\imath \lambda)} d\lambda.$$
Define the global holomorphic function
$$h_3(t):=- e^{-At} \int_{1-A}^0 e^{-t\lambda}
\frac{\dot{\Gamma}(A+\lambda)}{\Gamma(A+\lambda)}   d\lambda.$$
Since the singularities of $\Gamma(A-\imath \lambda)$ are $-iA-ik$,
$k=0,1,\dots$, we can
change the path of integration obtaining
$$I_1(t)=h_3(t)-\imath e^{-t}\int_0^\infty e^{\imath t\lambda}
\frac{\dot{\Gamma}(1-\imath \lambda)}{\Gamma(1-\imath \lambda)} d\lambda.$$
Assume $Re(t)>0,Im(t)>0$ for a moment.
Define
$$F(t):=\frac{1}{t}\int_0^\infty
\frac{\lambda}{e^\lambda-1}\frac{d\lambda}{\lambda+t}$$
and let $C$ be the Euler constant. It was shown by Cram\'er \cite{cramer19},
eq. (12), that
$$\imath \int_0^\infty e^{\imath t\lambda} \frac{\dot{\Gamma}(1-\imath
\lambda)}{\Gamma(1-\imath \lambda)} d\lambda=t\int_0^\infty e^{\imath t
\lambda} ln(\Gamma(1-\imath \lambda)) d(\imath\lambda)=\frac{C}{t}-F(t).$$
Moreover, there is a global meromorphic function
$$M(t):=F(t)-\frac{ln\:t}{e^{-t}-1},$$
where the branch of the logarithm is chosen such that $ln\:1=0$ and the cut is
at $(-\infty,0]$.
The poles of $M(t)$ are at $2\pi\imath \Z$ with residue $sign(k) \pi/2$ at
$t=k2\pi\imath$, $k\in\Z,k\not=0$.
Thus, we can write
$$I_1(t)=h_3(t)-e^{-t}\frac{C}{t}+ e^{-t} M(t)- \frac{ln\:t}{e^{t}-1}.$$
$I_1$ is a holomorphic function in $\C\setminus (-\infty,0]$.
The jump at the real negative axis is given by $\frac{2\pi\imath}{e^{t}-1}$,
if one crosses from above.

We evaluate
$$I_2(t):=\imath e^{-At}\int_0^\infty e^{\imath t\lambda}
\frac{\dot{\Gamma}(1/2+A-\imath \lambda)}{\Gamma(1/2+A-\imath \lambda)}
d\lambda.$$
Define the global holomorphic function
$$h_4(t):=e^{-At} \int_{1/2-A}^0 e^{-t\lambda}
\frac{\dot{\Gamma}(A+1/2+\lambda)}{\Gamma(A+1/2+\lambda)}   d\lambda.$$
Since the singularities of $\Gamma(A+1/2-\imath \lambda)$ are $-i(A+1/2)-ik$,
$k=0,1,\dots$, we can
change the path of integration obtaining
$$I_2(t)=h_4(t)+\imath e^{-t/2}\int_0^\infty e^{\imath t\lambda}
\frac{\dot{\Gamma}(1-\imath \lambda)}{\Gamma(1-\imath \lambda)} d\lambda.$$
Thus, we can write
$$I_2(t)=h_4(t)+e^{-t/2}\frac{C}{t}- e^{-t/2} M(t)+
\frac{e^{t/2}ln\:t}{e^{t}-1}.$$
$I_2$ is a holomorphic function in $\C\setminus (-\infty,0]$.
The jump at the real negative axis is given by $\frac{-2\pi\imath
e^{t/2}}{e^{t}-1}$,
if one crosses from above.

Now we consider
$$
e^{-At}\int_0^\infty e^{\imath t\lambda} \frac{\dot{G}(\imath
A+\lambda)}{G(\imath A+\lambda)} d\lambda.
$$
We evaluate
$$I_3(t):=\imath e^{-At}\int_0^\infty e^{\imath t\lambda}
\frac{\dot{\Gamma}(-A+\imath \lambda)}{\Gamma(-A+\imath \lambda)} d\lambda.$$
Assume, that $A$ is not an integer or half-integer.
Define the global holomorphic function
$$h_5(t):=e^{-At} \int^{A+1}_0 e^{t\lambda}
\frac{\dot{\Gamma}(-A+\lambda)}{\Gamma(-A+\lambda)}   d\lambda,$$
where the path avoids singularities   on the positive  imaginary part side.
The singularities of $\Gamma(-A+\imath \lambda)$ are $-iA+ik$, $k=0,1,\dots$.
We can
change the path of integration obtaining
$$I_3(t)=h_5(t)+\imath e^{t}\int_0^\infty e^{\imath t\lambda}
\frac{\dot{\Gamma}(1+\imath \lambda)}{\Gamma(1+\imath \lambda)} d\lambda.$$
Assuming $Re(t)<0$ for a moment
we can write
\begin{eqnarray*}\imath e^{t}\int_0^\infty e^{\imath t\lambda}
\frac{\dot{\Gamma}(1+\imath \lambda)}{\Gamma(1+\imath \lambda)}
d\lambda&=&e^{t}\int_0^\infty e^{t\lambda}
\frac{\dot{\Gamma}(1+\lambda)}{\Gamma(1+\lambda)} d\lambda\\
&=&-t e^{t}\int_0^\infty e^{t\lambda} ln\:\Gamma(1+\lambda) d\lambda\\
&=&t e^{t}\int_0^{-\infty} e^{-t\lambda} ln\:\Gamma(1-\lambda) d\lambda\\
&=&e^{t}(\frac{C}{t}+ F(t)).
\end{eqnarray*}
Now we can again apply a formula of Cram\'er (p. 114, top 3) and obtain
$$I_3(t)=h_5(t)+e^{t}\frac{C}{t}+e^{t} M(-t)-\frac{ln(-t)}{e^{-t}-1}.$$
$I_3$ is a holomorphic function in $\C\setminus [0,\infty)$.
The jump at the real positive axis is given by $\frac{2\pi\imath}{e^{-t}-1}$,
if one crosses from below.

We evaluate
$$I_4(t):=-\imath e^{-A t}\int_0^\infty e^{\imath t\lambda}
\frac{\dot{\Gamma}(-(A-1/2)+\imath \lambda)}{\Gamma(-(A-1/2)+\imath \lambda)}
d\lambda.$$
Define the global holomorphic function
$$h_6(t):=-e^{-A t} \int^{A+1/2}_0 e^{t\lambda}
\frac{\dot{\Gamma}(-(A-1/2)+\lambda)}{\Gamma((A-1/2)+\lambda)}   d\lambda,$$
where the path avoids the singularities  on the positive imaginary part side.
The singularities of $\Gamma(-(A-1/2)+\imath \lambda)$ are $-i(A-1/2)+ik$,
$k=0,1,\dots$. We can
change the path of integration obtaining
$$I_4(t)=h_6(t)-\imath e^{t/2}\int_0^\infty e^{\imath t\lambda}
\frac{\dot{\Gamma}(1+\imath \lambda)}{\Gamma(1+\imath \lambda)} d\lambda.$$
 Now we can again apply the formula of Cram\'er and obtain
$$I_4(t)=h_6(t)-e^{t/2}\frac{C}{t}-e^{t/2}
M(-t)+\frac{e^{-t/2}ln(-t)}{e^{-t}-1}.$$
$I_4$ is a holomorphic function in $\C\setminus [0,\infty)$.
The jump at the real positive axis is given by $\frac{-2\pi\imath
e^{-t/2}}{e^{-t}-1}$,
if one crosses from below.

Now we collect everything together.
Define the global holomorphic function
$$h(t):=h_1(t)+h_2(t)+rh_3(t)+rh_4(t)+rh_5(t)+rh_6(t)+2\pi\imath
\sum_{\sigma\in\Sigma, Re(\sigma)=0} m_\sigma (e^{-\imath t \sigma}-e^{\imath t
\sigma}).$$
\begin{theorem}
If $Im(t)>0$ then
\begin{eqnarray*}
2\pi\imath V(t)&:=&2\pi\imath (\tilde{\theta}_s(-t)-\theta_s(t))  \\
&=&h(t)+W(t)                 \\
&+&r(-e^{-t}\frac{C}{t}+e^{-t/2}\frac{C}{t/2}+
e^{t}\frac{C}{t}-e^{t/2}\frac{C}{t}\\
&&+e^{-t} M(t)- e^{-t/2} M(t)+e^{t}
 M(-t)-e^{t/2} M(-t)\\
&&-\frac{ln\:t}{e^{t}-1}+
\frac{e^{t/2}ln\:t}{e^{t}-1}-
\frac{ln(-t)}{e^{-t}-1}+
\frac{e^{-t/2}ln(-t)}{e^{-t}-1})\ .
\end{eqnarray*}
\end{theorem}
\begin{prop}
$$h(t)-h(-t)=0$$
\end{prop}
{\it Proof :}
We start with rewriting the integral defining $h_1$ using the functional
equation for $L$.
We have
\begin{eqnarray}
h_1(t)&=&-\imath \int_{-A}^A e^{-t(\lambda-\imath
0)}\frac{\dot{L}(\imath(\lambda-\imath 0))}{L(\imath(\lambda-\imath 0))}
d\lambda\nonumber  \\
      &=&-\imath \int_0^A e^{-t(\lambda-\imath
0)}\frac{\dot{L}(\imath(\lambda-\imath 0))}{L(\imath(\lambda-\imath 0))}
d\lambda \nonumber \\
      &&-\imath \int_{-A}^0 e^{-t(\lambda-\imath
0)}\frac{\dot{L}(\imath(\lambda-\imath 0))}{L(\imath(\lambda-\imath 0))}
d\lambda\label{jh1} ,
\end{eqnarray}
where we have indicated how the path avoids the singularities.
We rewrite (\ref{jh1}) as
$$-\imath \int_0^A e^{t(\lambda+\imath 0)}\frac{\dot{L}(-\imath(\lambda+\imath
0))}{L(-\imath(\lambda+\imath 0))} d\lambda$$
and apply the functional equation
$$\frac{\dot{L}(-\imath(\lambda+\imath 0))}{L(_\imath(\lambda+\imath
0))}=\frac{\dot{L}(\imath(\lambda+\imath 0))}{L(\imath(\lambda+\imath
0))}+\frac{\dot{G}(\imath(\lambda+\imath 0))}{G(\imath(\lambda+\imath 0))}-
  \frac{\dot{G}(-\imath(\lambda+\imath 0))}{G(-\imath(\lambda+\imath 0))}$$
in order to obtain
$h_1=h_7+rh_8$ with
\begin{eqnarray}
h_7(t)&=& -\imath \int_0^A \left( e^{-t(\lambda-\imath
0)}\frac{\dot{L}(\imath(\lambda-\imath 0))}{L(\imath(\lambda-\imath 0))}+
e^{t(\lambda+\imath 0)}\frac{\dot{L}(\imath(\lambda+\imath
0))}{L(\imath(\lambda+\imath 0))} \right)d\lambda
\nonumber\\
rh_8(t)&=& -\imath \int_0^A  e^{t(\lambda+\imath 0)}\left(
\frac{\dot{G}(\imath\lambda-0)}{G(\imath\lambda-
0)}+\frac{\dot{G}(-\imath\lambda+0)}{G(-\imath\lambda+0)} \right)d\lambda
.\label{jh2}
\end{eqnarray}
The anti-symmetrization $h_7(t)-h_7(-t)$ can be expressed in terms of the
singularities of $L$ on the interval $\imath (0,A)$. These singularities
are contributed by the scattering matrix $S$ as well as by the factor $G$.
We obtain
\begin{eqnarray*}h_7(t)-h_7(-t)&=&-2\pi\imath \sum_{\sigma\in\Sigma,
Re(\sigma)=0} 2m_\sigma (e^{-\imath t\sigma}-e^{\imath t\sigma})\\
&&-2\pi\imath r \sum_{0<k<A}(e^{-tk}-e^{tk})+2\pi\imath r
\sum_{0<k<A-1/2}(e^{-t(k+1/2)}-e^{t(k+1/2)}).
\end{eqnarray*}
Inserting the definition of $G$ into (\ref{jh2}) we obtain
\begin{eqnarray*}
h_8(t)&=&\int_0^A e^{t(\lambda+\imath 0)} (\frac{\dot{\Gamma}(-\lambda-\imath
0)}{\Gamma(-\lambda-\imath 0)} - \frac{\dot{\Gamma}(1/2 - \lambda-\imath
0)}{\Gamma(1/2 -\lambda-\imath 0)}  \\
      &&-\frac{\dot{\Gamma}(\lambda+\imath 0)}{\Gamma(\lambda+\imath 0)}   +
\frac{\dot{\Gamma}(1/2+\lambda+\imath 0)}{\Gamma(1/2+ \lambda+\imath 0)}).
\end{eqnarray*}
We also have
\begin{eqnarray*}
h_3(t)&=&-\int_1^A e^{-t\lambda} \frac{\dot{\Gamma}(\lambda)}{\Gamma(\lambda)}
d\lambda\\
h_4(t)&=&\int_{1/2}^A e^{-t\lambda}
\frac{\dot{\Gamma}(1/2+\lambda)}{\Gamma(1/2+\lambda)}  d\lambda\\
h_5(t)&=&\int_{-1}^A e^{-t\lambda} \frac{\dot{\Gamma}(-\lambda+\imath
0)}{\Gamma(-\lambda+\imath 0)}  d\lambda\\
h_6(t)&=&-\int_{-1/2}^A e^{-t\lambda} \frac{\dot{\Gamma}(1/2-\lambda+\imath
0)}{\Gamma(1/2-\lambda+\imath 0)}  d\lambda.
\end{eqnarray*}
Thus, defining
$$h_9(t):=h_8(t)-h_3(-t)-h_4(-t)-h_5(-t)-h_6(-t)$$
and anti-symmetrizing we obtain
$$h_9(t)-h_9(-t)=2\pi\imath \sum_{0<k<A}(e^{-tk}-e^{tk})-2\pi\imath
\sum_{0<k<A-1/2}(e^{-t(k+1/2)}-e^{t(k+1/2)}).$$
Since $h_2$ is symmetric,
$$h(t)-h(-t)=h_7(t)-h_7(-t)+r(h_9(t)-h_9(-t))+2\pi\imath \sum_{\sigma\in\Sigma,
Re(\sigma)=0} 2m_\sigma (e^{-\imath t \sigma}-e^{\imath t \sigma})=0$$
and the proposition follows.
$\Box$\newline

We can now write down a distributional trace formula on $(-\infty,0)$
containing
$\theta(t+
\imath 0)+\tilde{\theta}(t-\imath 0)$.
In fact, for $t>0$
\begin{eqnarray*}
2\pi\imath(\theta_s(
-t+\imath 0)+\tilde{\theta}_s(-t-\imath
0))&=&2\pi\imath(\theta_s(t+\imath 0)+
\tilde{\theta}_s(t-\imath 0)+V(t+\imath
0)-V(-t+\imath 0))\\
&=&2\pi
\imath(\theta_s(t+\imath 0)+
\tilde{\theta}_s(t-\imath 0))\\
&+&W(t+\imath 0)-W(-t+\imath
0)\\
&+&2\pi\imath r
(\frac{1}{e^{-t}-1}-\frac{e^{-t/2}}{e^{-t}-1})\ .
\end{eqnarray*}
By (\ref{er3}) we have
$$W(t+\imath 0)-W(-t+\imath 0)=-2\pi\imath  \sum_{p\in R} ln|p| c(p)
\delta(t-ln|p|).$$
We combine these two equations with the distributional trace formula and obtain
\begin{theorem}
The following equation of distributions on $(-\infty,0)$ holds :
\begin{eqnarray*}
\theta(t+\imath 0)+\tilde{\theta}(t-\imath 0)&=&-r \frac{e^{t/2}}{e^{t}-1}\\
                   &-& \frac{vol(M)}{4\pi}\frac{cosh(t/2)}{sinh^2(t/2)}\\
                                             &+&\sum_{c\mbox{\scriptsize
-closed geodesic}}\frac{l_c}{2n_csinh(l_c/2)}\delta(t+l_c)\\
                                             &-&\sum_{p\in R} ln|p| c(p)
\delta(t+ln|p|) .
\end{eqnarray*}
\end{theorem}
We can now apply exactly the same technique as for Theorem \ref{th3} in order
to prove
\begin{theorem}
$\theta(t)$ has a meromorphic continuation from the lower half-plane
to another sheet of the upper half-plane across $(-\infty,0)$.
For $Re(t)>0$ this extension $\theta_1$ is given by
\begin{eqnarray*}
\theta_1(t)&=&\theta(t)+r(\frac{1}{1-e^{-t}} - \frac{e^{t/2}}{e^{t}-1})\\
           &=&\theta(t)+r\frac{e^{t/2}-1}{2 sinh(t/2)}.
\end{eqnarray*}
The singularities on the negative real axis are first order poles
in the points $-l_c$, $c$ a closed geodesic, with residue $\frac{l_c}{4\pi
n_csinh(l_c/2)}$
and in the points $-ln|p|$, $p\in R$, with residue $\frac{-ln|p| c(p)}{2\pi}$
(if two such points
coincide, the residues add up).
Moreover, it has first order poles in the points $(4k-2)\pi\imath $,
$k=1,2,\dots$,
with residue $-2 r$.
\end{theorem}

\section{Conclusion}
We have seen that $\theta(t)$ extends meromorphically to the Riemann
surface of the logarithm. The difference between two sheets is
$$\theta_1(t)-\theta(t)=r\frac{e^{t/2}-1}{2 sinh(t/2)}\ . $$
Thus, the modified theta-function
$$\Theta(t):=\theta(t)+
\frac{r}{2\pi\imath}ln(t) \frac{e^{t/2}-1}{2 sinh(t/2)} $$
admits a meromorphic extension to the complex plane.
\begin{theorem}
The meromorphic extension of the modified theta function has the following
singularities:
\begin{itemize}
\item First order poles at $l_c$, $c$ a closed geodesic, with residue
$\frac{l_c}{4\pi n_csinh(l_c/2)}$,
\item first order poles at $l_c$, $c$ a closed geodesic, with residue
$\frac{l_c}{4\pi n_csinh(l_c/2)}$
      and at $-ln|p|$, $p\in R$, with residue $\frac{-ln|p| c(p)}{2\pi}$ (if
two such points
      coincide, the residues add up)
\item first order poles at $(4k-2)\pi\imath $, $k=1,2,\dots$, with residue
$-r/2-ln(\pi(4k-2))\frac{r}{\pi\imath}$
\item second order poles at $(4k+2)\pi\imath $, $k=-1,-2,\dots$, with residue
$-3r/2-ln(\pi(4k+2))\frac{r}{\pi\imath} $ and second Laurent coefficient
$-\frac{vol(M)}{\pi}$
\item and second order poles at $4k \pi \imath $, $k=-1,-2,\dots$, with residue
$r$ and second Laurent coefficient $-\frac{vol(M)}{\pi}$.
\end{itemize}
\end{theorem}
Clearly, $\Theta$ has at most a pole of second order at $t=0$.

\bibliographystyle{plain}

\end{document}